\pdfoutput=1
\documentclass[11pt,oneside]{amsart}
\usepackage[leqno]{amsmath}
\usepackage{amssymb,amsthm} 
\usepackage[parfill]{parskip} 
\usepackage{libertine} %loads biolinum ss, similar to optima
\usepackage[T1]{fontenc}
\usepackage{textcomp}
\usepackage[varqu,varl]{inconsolata} %tt
\usepackage{amsmath,amsthm}
\usepackage[libertine,bigdelims,vvarbb]{newtxmath}
\usepackage[supstfm=libertinesups,%
  supscaled=1.2,%
  raised=-.13em]{superiors}
\usepackage{graphicx} 
\usepackage{calc}
\newcommand{\beq}{\begin{equation}}
\newcommand{\eeq}{\end{equation}}
\numberwithin{equation}{section} %restart every new section
\newtheorem{thm}[equation]{Theorem}
\newtheorem{lem}[equation]{Lemma}

\theoremstyle{remark}

\theoremstyle{definition}
\newtheorem{defn}[equation]{Definition}

\title{ZFC independence and subset sum}
\usepackage[paperwidth=7in, paperheight=10in, textwidth=4.75in, textheight=8in, includehead=true]{geometry}
\usepackage{url}
\usepackage[pdfborderstyle={/S/U/W 1}]{hyperref}

\def\setmax{\operatorname{setmax}}

\let\emptyset\emptysetAlt
\hypersetup{
pdftitle={ZFC},
pdfauthor={\textcopyright\ S. G. Williamson}
}
\author{S. Gill Williamson}
\thanks{Department of Computer Science and Engineering, 
University of California San Diego; \url{http://cseweb.ucsd.edu/~gill/}.
{\bf Keywords:} lattice exit models,  ZFC independence, 
order type equivalence, regressive regularity, subset sum instances, P=NP
}
\date{08/25/2017}                                           
\begin{document}
\begin{abstract}
Let $Z$ be the integers and $N$ the nonnegative integers.
and let $G=(N^k,\Theta)$ be a ``$\max$-downward'' digraph.  
We study sets of functions 
$H=\{h_D^\rho\mid h_D^\rho:D\rightarrow N, D\subset N^k ,|D|<\infty\}$ parameterized by sets 
$\rho=\{\rho_D\mid\rho_D:D \rightarrow N, \rho_D(x)\geq \min(x)\}$  and defined recursively using the structure of $G$. 
We prove the sets $H$ satisfy: for $p\geq 2$ there exists 
$\hat{E}\subset N$, $|\hat{E}|=p$,  $\hat{E}^k\subseteq \hat{D}$, $h_{\hat{D}}^\rho$  ``regressively regular'' over ${\hat{E}}$ (def~\ref{def:regreg}). 
We show that this theorem is independent of ZFC (thm~\ref{thm:jfhcaplog}).
We define a parameterized family of finite subsets of $Z$, 
$S^\rho_{F,G}(k,t,E,p,D),\;$ (def~\ref{def:setsdis}).
Fixing all parameters but $E$, $p$ and $D$ and using regressive regularity to uniquely associate $\hat{D}$ and $\hat{E}$  with $p$ we obtain 
sequences of sets of integers $\{S(\hat{E}, \hat{D}): p=2,3,\ldots \}$.  We show any such set of instances to the subset sum problem can be solved in time $O(p^t)$ for some $t\geq 1$ (thm~\ref{cor:subsumpoly}).

\end{abstract}

\maketitle
%\section{A}
%\subsection{}
%\section*{}\hbox to\textwidth{\hfil \textbf{Abstract}\hfil}{\small this is my abstract}

\section{Introduction}
Basic references are Friedman~\cite{hf:alc} and Williamson\cite{gw:lem}.
We extend a result derived in \cite{gw:lem} in order to make a connection between {\em regressive regularity} type ZFC independence and the subset sum problem.  
In particular, we extend definition~4.6 and theorem~4.8 of \cite{gw:lem} to
definition~\ref{def:cap} ($D$  capped by $E^k\subset D$) and theorem~\ref{thm:jfhcap} (regressive regularity of 
$h_D^\rho$, capped version).
The latter  results are extended to definition~\ref{log:bound} ($\rho_D$ log bounded over $E$) and theorem~\ref{thm:jfhcaplog} (regressive regularity of $h_D^\rho$, log bounded version). Both theorem \ref{thm:jfhcap} theorem~\ref{thm:jfhcaplog} are independent of ZFC (assuming consistency). 

In definition~\ref{def:setsdis} we define an uncountably infinite family, 
$S^\rho_{F,G}(k,t,E,p,D)$, of finite subsets of integers,   We call these finite subsets of integers
{\em sets of displacements}. 
Each set of displacements is of the form $\delta_E h_D^\rho E^k_l \cup \delta_E h_D^\rho \mathrm{diag}(E^k)$ (notation to be explained below).
These sets of displacements are constructed to closely reflect the structure of theorem~\ref{thm:jfhcaplog}.
Putting off exact definitions for the moment, the parameters of 
$S^\rho_{F,G}(k,t,E,p,D)$ are as follows:  $G=(N^k,\Theta)$, $\;k\geq 2,\;$ ranges over all downward directed lattice graphs of dimension $k$, $t\geq 1$, $F:N^k\times (N^k\times N)^r\rightarrow N,\;$ $r\geq 1,\;$ ranges over all ``partial selection functions,''
$E\subset N$ ranges over all finite subsets, $|E|=p\geq 2$, $D\subset N^k$ ranges over all  $D$ ``capped by'' $E^k.$
%and $h_D^\rho$ has $\rho_D$ (indicated by the exponent $\rho$) $t$-log bounded ($\rho_D\in \rm{LOG}(k,E,p,D,t)$).

Our main result, theorem~\ref{cor:subsumpoly}, states that 

``For fixed $F, G, k,$ consider sets of instances 
$\{\delta_E h_D^\rho E^k_l \cup \delta_E h_D^\rho \mathrm{diag}(E^k): E, p, D\}.$
For each $p$ there exists $\hat{E}$ and $\hat{D}$ such that
the subset sum problem for
$$\{\delta_{\hat{E}} h_{\hat{D}}^\rho {\hat{E}}^k_l \cup \delta_{\hat{E}} h_{\hat{D}}^\rho \mathrm{diag}(\hat{E}^k): p=2, 3, \ldots\}$$ is solvable in time $O(p^t)$ for some $t$.''

Here $p=|\hat{E}|$ is a measure of the size of the instance. 

Our only proof of theorem~\ref{cor:subsumpoly} is by using the ZFC independent 
theorem~\ref{thm:jfhcaplog}.
If ``subset sum solvable in polynomial time'' could be proved in ZFC then 
a ZFC proof of 
theorem~\ref{cor:subsumpoly} would follow (no need for theorem~\ref{thm:jfhcaplog}).  We conjecture that 
 theorem~\ref{cor:subsumpoly} is itself independent of ZFC.  
 If so, ``subset sum solvable in polynomial time'' would be independent of ZFC.  
  
\section{Elementary background}
Let $N$ be the set of nonnegative integers and $k\geq 2$.
For $z=(n_1, \ldots, n_k)\in N^k$, $\max\{n_i\mid i=1,\ldots, k\}$ will be
denoted by $\max(z)$.  Define $\min(z)$ similarly. 
\begin{defn}[Downward directed graph]
\label{def:dwndrctgrph}
Let $G=(N^k,\Theta)$ (vertex set $N^k$, edge set $\Theta$) be a directed graph.
If every $(x,y)$ of $\Theta$ satisfies $\max(x) >\max(y)$  then we call $G$ a {\em downward directed lattice graph}. 
For $z\in N^k$, let $G^z = \{x: (z,x)\in \Theta\}$ denote the vertices of $G$ {\em adjacent} to $z$.
\end{defn}
{\bf All lattice graphs that we consider will be {\em downward directed}.}\\
%FIGURE DOWNWARD DIRECTED 
%\begin{figure}[h]
%\begin{center}
%\includegraphics[scale=.85]{./Figures/Downdir}
%\caption{Downward directed, $k=2$}
%\label{fig:downdir}
%\end{center}
%\end{figure}
%END FIGURE
%DEFINITION VERTEX INDUCED SUBGRAPH
\begin{defn} [\bf Vertex induced subgraph $G_D$]
For $D\subset N^k$  let 
$G_D = (D, \Theta_D)$ be the subgraph of $G$ with vertex set $D$ and edge set 
$\Theta_D =\{(x,y)\mid (x,y)\in \Theta,\, x, y \in D\}$. We call $G_D$ the {\em subgraph of $G$ induced by $D$}. 
\label{def:vertexinduced}
\end{defn}
%END DEFINITION
\vskip .25in  
%\begin{defn} [\bf Path and terminal path in $G_D$]
%\label{def:terminalpath}
%For $t\geq 2$, a sequence of distinct vertices of $G_D$, $(x_1, x_2, \ldots, x_t)$, is a {\em path}
%of length $t$ in $G_D$  if $(x_i, x_{i+1}) \in \Theta_D,\, i=1, \ldots\,, t-1$.  For $
%x\in D$, $(x)$ as a path of length $1$. The path $(x_1, x_2, \ldots, x_t)$ is {\em terminal}
%if there is no path of the form  $(x_1, x_2, \ldots, x_t, x_{t+1})$. 
%We say $x$ is a {\em terminal vertex} of $G_D$ if the path $(x)$ is a terminal path in $G_D$.
%\label{def:downpath}
%\end{defn}

\begin{defn}[\bf Cubes and Cartesian powers in $N^k$]
The set  $E_1\times\cdots\times E_k$, where $E_i\subset N$, $|E_i|=p$, $i=1,\ldots, k,$   are $k$-cubes of length $p$.  If $E_i = E, i=1,\ldots, k,$ then this cube is  $E^k =\times^k E$, the $k$th Cartesian power of $E$.
\label{def:cubespowers}
\end{defn}

\begin{defn}[\bfseries Equivalent ordered $k$-tuples]
\label{def:ordtypeqv}
Two k-tuples in $N^k$, $x=(n_1,\ldots,n_k)$ and $y=(m_1,\ldots,m_k)$, are  
{\em order equivalent tuples $(x\, ot \,y)$} if 
$\{(i,j)\mid n_i < n_j\} =  \{(i,j)\mid m_i < m_j\}$ and  $\{(i,j)\mid n_i = n_j\} =  \{(i,j)\mid m_i =m_j\}.$  
\end{defn}
Note that $ot$ is  an equivalence relation on $N^k$.
The standard SDR (system of distinct representatives) for the $ot$  
equivalence relation is gotten by replacing $x=(n_1,\ldots,n_k)$ by 
$\mathbf{r}(x)= (\mathbf{r}_{S_x}(n_1),\ldots,\mathbf{r}_{S_x}(n_k))$ where
 $\mathbf{r}_{S_x}(n_j)$ is the rank of $n_j$ in $S_x = \{n_1, \ldots, n_k\}$ (e.g,
$x=(3, 8, 5, 3, 8)$, $S_x = \{3, 5, 8\}$, $\mathbf{r}(x)=(0, 2, 1, 0, 2)$).
The number of equivalence classes is $\sum_{j=1}^k \sigma(k, j)\leq k^k$ where
$\sigma(k,j)$ is the number of surjections from a $k$ set to  a $j$ set.
We use ``$x\,ot\,y$'' and ``$x,\,y$ of order type $ot$'' to mean $x$ and $y$ belong to the same order type equivalence class.

%\begin{defn}[$\Ki$ notation]
%Define $\Ki(statement) = 0\;{\rm if}\,statement\, {\rm false}$,
%and $\Ki(statement) = 1\;{\rm if}\,statement\, {\rm true}$.
%\end{defn}
\section{Basic definitions and theorems}

\begin{defn}[\bf regressive value]
Let $X\subseteq N^k$ and $f:X\rightarrow Y\subseteq N$.  An integer $n$
is a  {\em regressive value} of $f$ on $X$
if there exist $x$ such that $f(x)=n<\min(x)$ .
\end{defn}

\begin{defn}[\bf field of a function and reflexive functions]
For $A\subseteq N^k$ define ${\rm field}(A)$ to be the set of all coordinates of elements of $A$.  A function $f$ is reflexive in $N^k$ if 
${\rm domain}(f) \subseteq N^k$ and  ${\rm range}(f) \subseteq {\rm field}({\rm domain}(f))$.
\end{defn}

\begin{defn}[the set of functions $T(k)$ ]
$T(k)$ denotes all reflexive functions with finite domain: $|{\rm domain}(f)|<\infty$.
\end{defn}

\begin{defn} [\bf full and jump free] 

Let $Q\subset T(k).$

\begin{enumerate}
\item {\bf full:}  $Q$ is a {\em full} family of functions on $N^k$ if for every finite subset 
$D\subset N^k$ there is at least one function $f$ in $Q$ whose domain is $D$.

\item{\bf jump free:} For $D\subset N^k$ and $x\in D$ define $D_x = \{z\mid z\in D,\, \max(z) < \max(x)\}$. 
Suppose that for all $f_A$ and $f_B$  in $Q$, where $f_A$ has domain $A$ and $f_B$ has domain $B$,  the conditions
 $x\in A\cap B$, $A_x \subseteq B_x$, and $f_A(y) = f_B(y)$ for all $y\in A_x$ imply that 
$f_A(x) \geq f_B(x)$.  Then $Q$ will be called a {\em jump free} family of functions on $N^k$.  
\end{enumerate}
\label{def:fullrefjf}
\end{defn} 

\begin{defn}[\bf Regressively regular over $E$]
\label{def:regreg}
Let $k\geq 2$, $D\subset N^k$, $D$ finite, $f: D\rightarrow N$. 
We say $f$ is {\em regressively regular} over 
$E$, $E^k\subset D$, if for each order type equivalence class $ot\,$ of $k$-tuples of $E^k$ either (1) or (2) occurs:
\begin{enumerate}
\item{\bf constant less than min $\mathbf{E}$:}  For all $x, y\,\in E^k$ of order type $ot$, $f(x)=f(y)< \min(E)$ 
\item{\bf greater than min:}   For all $x\in E^k$ of order type $ot$ $f(x)\geq \min(x).$
\end{enumerate}
\end{defn}   

\begin{figure}[h]
\begin{center}
\includegraphics[scale=.85]{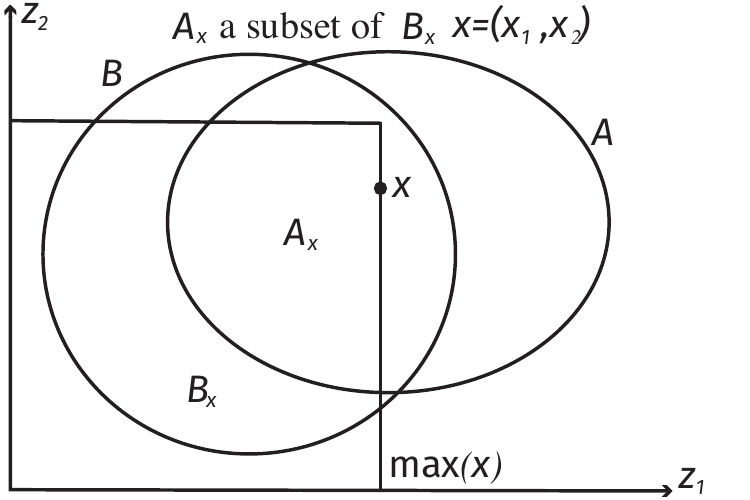}
\caption{Basic jump free condition~\ref{def:fullrefjf}}
\label{fig:jfvenn}
\end{center}
\end{figure}
\begin{thm}[\bf Jump free theorem (\cite{hf:alc}, \cite{hf:nlc})] 
\label{thm:jumpfree}
Let $p, k\geq 2$ and $S\subseteq T(k)$ be  a full and jump free family.
Then some $f\in S$ has at most $k^k$ regressive values on some 
$E^k \subseteq {\rm domain}(f)$, $|E| = p$.  
In fact, some $f\in S$ is regressively regular over some $E$ of cardinality $p$.
\end{thm}

We use ZFC for the axioms of set theory, Zermelo-Frankel  plus the axiom of choice (see Wikipedia).  
 The jump free theorem can be proved in 
 ZFC + ($\forall n$)($\exists$ $n$-subtle cardinal)
 but not in ($\exists$ $n$-subtle cardinal) for any fixed $n$ (assuming this theory is consistent).
A proof is in Section 2 of \cite{hf:alc},
``Applications of Large Cardinals to Graph Theory,'' October 23, 1997, No. 11 of  
\href{https://u.osu.edu/friedman.8/foundational-adventures/downloadable-manuscripts/}{Downloadable Manuscripts}.

\section{Large scale regularities}

\begin{defn} [{\bf Partial selection}]
\label{def:partselect}
A function $F$ with domain a subset of $X$ and range a subset of $Y$ will be called a {\em partial function}
from $X$ to $Y$ (denoted by $F: X\rightarrow Y$).  If $z\in X$ but $z$ is not in the domain of $F$, we say 
$F$ is {\em not defined} at $z$.
Let $r
\geq 1$.  A partial function 
$F: N^k\times(N^k \times N)^r \rightarrow N$
will be called a {\em partial selection} function ~\cite{hf:alc} if whenever 
$F(x, ((y_1,n_1), (y_2,n_2), \ldots (y_r,n_r)))$ is defined we have 
$F(x, ((y_1,n_1), (y_2,n_2), \ldots (y_r,n_r))) = n_i$ for some $1\leq i \leq r$.
\end{defn}

\begin{defn} [\bf Committee model {\bf $\hat{s}_D$}\cite{hf:alc}, \cite{gw:lem}]
\label{def:chanlabel}
Let $r\geq 1,$ $k\geq 2$, $G=(N^k,\Theta),$ $G_D = (D, \Theta_D)$, $D$ finite, $G_D^z= \{x\mid (z,x)\in \Theta_D\}.$ 
Let $F: N^k\times(N^k \times N)^r \rightarrow N$ be a partial selection function.  
We define $\hat{s}_D(z)$ recursively (on $\max(z)$) on $D$ as follows. Let
\[
\Phi^D_z = \{ F[z, (y_1,n_1), (y_2,n_2), \ldots, (y_r,n_r)],\;y_i \in G^z_D\}
\]
be the set of defined values of $F$  where  
$n_i=\hat{s}_D(y_i)$ if $\Phi^D_{y_i}\neq\emptyset$ and
$\;n_i=\min(y_i)$ if $\Phi^D_{y_i}=\emptyset.\;$
If $\Phi^D_z=\emptyset$,  define $\hat{s}_D(z) = \max(z)$.
If $\Phi^D_z\neq\emptyset$,  define $\hat{s}_D(z)$ to be the minimum over $\Phi^D_z$.
\end{defn}

{\bf NOTE:} If $\Phi^D_z\neq \emptyset$ then an induction on $\max(z)$ shows 
 $\hat{s}_D(z) < \max(z).$ Recall that $(G,\Theta)$ is downward.
Thus, $\Phi^D_z=\emptyset$ iff $\hat{s}_D(z) = \max(z)$ 
(see lemma~\ref{lem:shatvsh}).

\begin{thm}[\bfseries Large scale regularities for $\hat {s}_D$]
\label{thm:jfhats}
Let $r\geq 1$, $p, k\geq 2$.
$S=\{\hat {s}_D\mid D\subset N^k,\;|D|<\infty\}$. Then some 
$f\in S$ has at most $k^k$ regressive values over some 
$E^k \subseteq {\rm domain}(f)$, $|E|=p.$ 
In fact, some $f\in S$ is regressively  regular over some $E$ of cardinality $p$.

\begin{proof}
Recall \ref{thm:jumpfree}. Let  $S=\{\hat {s}_D\mid D\subset N^k,\;|D|<\infty\}.$
$S$ is obviously full and reflexive. 
We show $S$ is jump free.
We show for all $\hat {s}_A$ and $\hat {s}_B$  in $S$,  the conditions
 $x\in A\cap B$, $A_x \subseteq B_x$, and $\hat {s}_A(y) =\hat {s}_B(y)$ for all $y\in A_x$ imply that 
$\hat {s}_A(x) \geq \hat {s}_B(x)$.  (i.e., $S$ is  {\em jump free}). 
If $\Phi^A_x = \emptyset$ then $\hat {s}_A(x)=\max(x)\geq \hat {s}_B(x).$
Assume $\Phi^A_x \neq \emptyset.$ 
Let $n=F[x, (y_1,n_1), (y_2,n_2), \ldots (y_r,n_r)]\in \Phi^A_x$ 
(note that $y_i\in G_A^x \subseteq G_B^x$) where 
$n_i=\hat {s}_A(y_i)$ if $\hat {s}_A(y_i)<\max(y_i)$ 
(i.e., $\Phi^A_{y_i} \neq \emptyset$ see {\bf NOTE} after definition~\ref{def:chanlabel}) and 
$n_i=\min(y_i)$ if $\hat {s}_A(y_i)=\max(y_i).$
But $\hat {s}_A(y_i) =\hat {s}_B(y_i)$, $i=1, \ldots, r,$ implies
$n\in \Phi^B_x$ and thus $\Phi^A_x \subseteq \Phi^B_x$ and
$
\hat {s}_A(x)=\min(\Phi^A_x) \geq \min(\Phi^B_x)=\hat {s}_B(x).
$
\end{proof}
\end{thm}

\begin{figure}[h]
\begin{center}
\includegraphics[scale=.7]{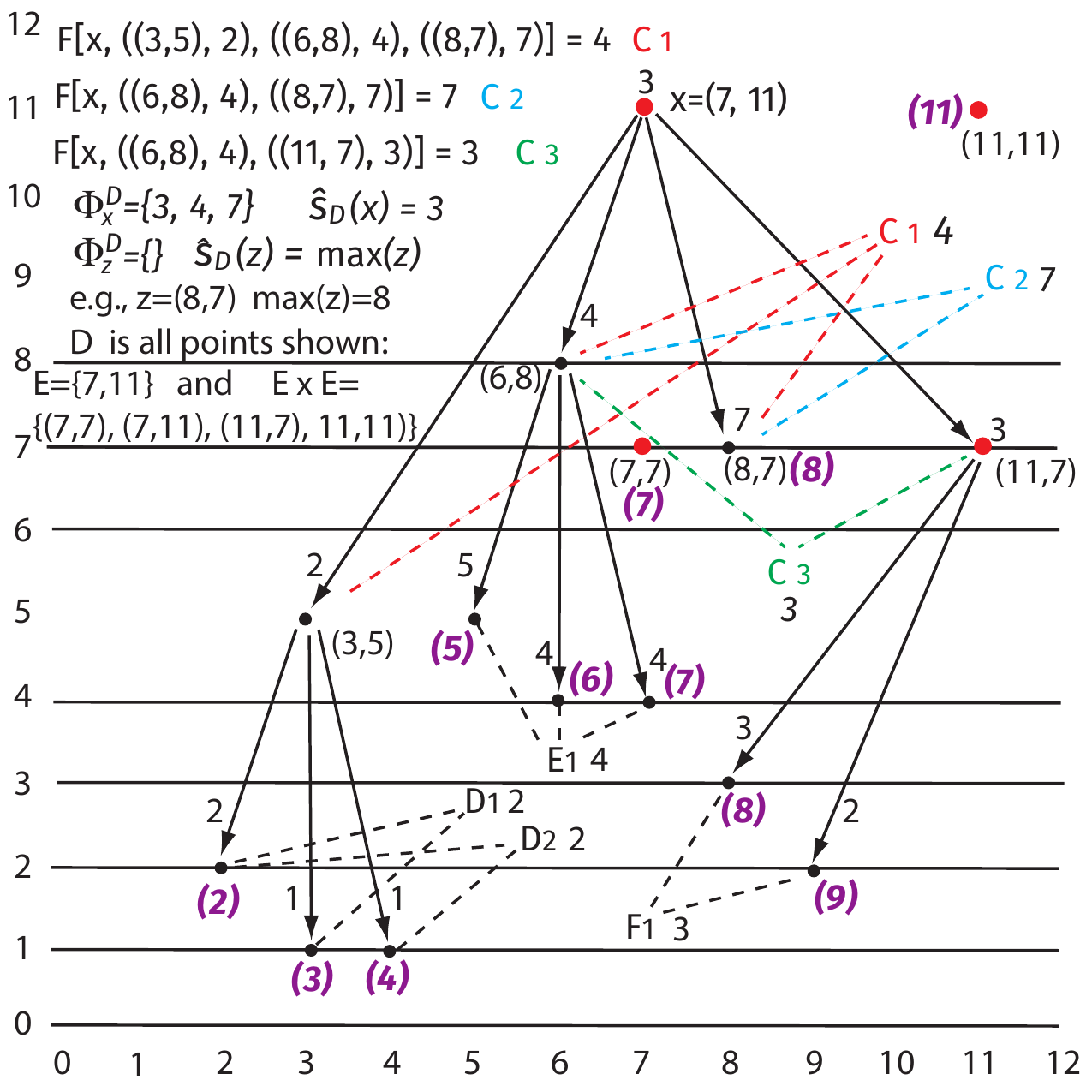}
\caption{An example of $\hat{s}_D$}
\label{fig:univs}
\end{center}
\end{figure}

As an example of computing $\hat{s}_D$, consider figure~\ref{fig:univs}.
The computation is recursive on the $\max$ norm.  
The values of the terminal vertices where  $\Phi^A_x = \emptyset$ are shown in parentheses,  left to right: (2), (3), (4), (5), (6), (8), (8), (9).  These numbers are  $\max((a,b))$ for each terminal vertex $(a,b)$.
Partial selection functions are of the form 
$F:N^2 \times (N^2\times N)^r \rightarrow N$ ($r=2, 3$ here).  
In particular we have
$F[x, ((3,5),2), ((6,8),4), ((8,7),7)]=4,\;\;$
$F[x, ((6,8),4), ((8,7),7)]= 7$, and\\
$F[x, ((6,8),4), ((11,7),3)]= 3$.
Intuitively, we think of these as (ordered) committees reporting values to the boss, $x=(7,11).$
The first committee, $\rm{C}1$, consists of subordinates, $(3,5),(6,8), (8,7)$
reporting respectively $2, 4, 7$.  
The committee decides to report $4$ (indicated by  $\rm{C}1\;4$ in  
figure~\ref{fig:univs}). 
The recursive construction starts with terminal vertices reporting their minimal coordinates.
But, the value reported by each committee is not, in general, the actual minimum of the reports of the individual members. 
Nevertheless, the boss, $x=(7,11)$,  {\em always} takes the minimum of the values reported by the committees.  
In this case the values reported by the committees are $4, 7, 3$ the boss takes $3$ (i.e., $\hat{s}_D(x)=3$ for the boss, $x=(7,11)$).
Note that a function like $F((7,11), ((6,8),4), ((8,7),7)$ where $r=2$, can be padded
to the case $r=3$\\
(e.g., $F((7,11), ((6,8),4), ((8,7),7), ((8,7),7))).$ 

Observe in figure~\ref{fig:univs} that the values in parentheses, (2), (3), (4), (5), (6), (8), (8), (9), don't figure into the recursive construction of $\hat{s}_D$.  They immediately pass their minimum values on to the computation: 2, 1, 1, 5, 4, 4, 7, 3.
This leads to the following generalization of definition~\ref{def:chanlabel}.

\begin{defn} [{\bf $h^{\rho}_D$ for $G_D$}]
\label{def:genhat}
Let $r\geq 1,$ $k\geq 2$, $G=(N^k,\Theta),$ $G_D = (D, \Theta_D)$, $D$ finite, $G_D^z= \{x\mid (z,x)\in \Theta_D\}.$ 
Let $F: N^k\times(N^k \times N)^r \rightarrow N$ be a partial selection function.
Let $\rho=\{\rho_{D}\mid \rho_D:D \rightarrow N,\, \min(x)\leq \rho_{D}(x),\,x\in D\}$.
We define $h_D^\rho(z)$ recursively (on $\max(z)$) on $D$ as follows. 
Let
\[
\Phi^D_z = \{ F[z, (y_1,n_1), (y_2,n_2), \ldots, (y_r,n_r)],\;y_i \in G^z_D\}
\]
be the set of defined values of $F$  where  
$n_i=h^{\rho}_D(y_i)$ if $\Phi^D_{y_i}\neq\emptyset,$ and
$\;n_i=\min(y_i)$ if $\Phi^D_{y_i}=\emptyset.\;$
If $\Phi^D_z=\emptyset$,  define $h^{\rho}_D(z) = \rho_{D}(z)$.
If $\Phi^D_z\neq\emptyset$,  define $h^{\rho}_D(z)$ to be the minimum over 
$\Phi^D_z$.
Note that $\rho_D$ need not be reflexive on $D$.
\end{defn}

\begin{lem}[{\bf Compare} $\hat{s}_D,\, h^{\rho}_D$]
\label{lem:shatvsh}
For all $z\in D$, $h^{\rho}_D(z)=\hat{s}_D(z)<\max(z)$ if $\Phi^D_z\neq \emptyset$ and
$h_D^\rho (z) = \rho_D(z)$, $\;\hat{s}_D(z)=\max(z)$  if $\Phi^D_z=\emptyset$.
Let $E$ be of cardinality $p\geq 2$.
Then $\hat{s}_D$ regressively  regular over $E$ iff $h^{\rho}_D$ regressively  regular over $E$. 
\begin{proof}
We use induction on $\max$. 
Let $D_a=\{x\mid x\in D,\, \max(x)=a\}$.
Let $m_0 < m_1 < \cdots < m_q$ be the integers $n$ such that $D_n\neq \emptyset$.
If $z\in D_{m_0}$ then the set of adjacent vertices $G_D^z = \emptyset$.
Thus, $\Phi^D_z = \emptyset$ and $h^{\rho}_D(z) = \rho_D(z)$, $\hat{s}_D(z)=\max(z)$.
Consider $z\in D_{m_1}$.  
If $\Phi^D_z = \emptyset$ then $h^{\rho}_D(z) = \rho_D(z)$, $\hat{s}_D(z)=\max(z)$.
Assume $\Phi^D_z \neq \emptyset$ and let $n=F[z, (y_1,n_1), (y_2,n_2), \ldots, (y_r,n_r)],\;y_i \in G^z_{D}.$  
But, $\Phi^D_{y_i} = \emptyset$ for all $y_i\in D_{m_0}$ implies 
$n_i = \min(y_i), i=1\ldots r$.  
This observation is the same whether computing  $h_D^\rho (z)$ or $\hat{s}_D(z)$.
Thus, for $z\in D_{m_1}$, $h_D^\rho (z) = \hat{s}_D(z)<\max(z)$ if $\Phi^D_z \neq \emptyset$ and
$h_D^\rho (z) = \rho_D(z)$, $\;\hat{s}_D(z)=\max(z)$ otherwise.

Assume, for $y\in D_{m_t}$ where $t<j$, $h_D^\rho (y) = \hat{s}_D(y)<\max(y)$ if $\Phi^D_y \neq \emptyset$
and
$h_D^\rho (y) = \rho_D(y)$, $\hat{s}_D(y)=\max(y)$ otherwise.
Let $z\in D_{m_{j}}$.
If $\Phi^D_z = \emptyset$ then $h^{\rho}_D(z) = \rho_D(z)$ and $\hat{s}_D(z)=\max(z)$.
Let  $n=F[z, (y_1,n_1), (y_2,n_2), \ldots, (y_r,n_r)]$ for $y_i \in G^z_{D}$  
and thus $\Phi^D_z \neq \emptyset$.
By induction, if $\Phi^D_{y_i} \neq \emptyset$ then $n_i = h_D^\rho(y_i) = \hat{s}_D(y_i)<
\max(y_i).$ 
If $\Phi^D_{y_i} = \emptyset$ then $n_i = \min(y_i)$ in computing either 
$h_D^\rho (z)$ or  $\hat{s}_D(z)$.
Thus,  $h_D^\rho (z) = \hat{s}_D(z)<\max(z)$ if $\Phi^D_z \neq \emptyset$
and
$h_D^\rho (z) = \rho_D(z)$, $\hat{s}_D(z)=\max(z)$ otherwise.
 
Finally, we consider regressive regularity.  {\bf First} we show for all $x, y\,\in E^k$ of order type $ot$, 
$\hat{s}_D(x)=\hat{s}_D(y)< \min(E)$ 
 if and only if $h^\rho_D(x)=h^\rho_D(y)< \min(E)$.  In ether case  $\Phi^D_x \neq \emptyset$
and $\Phi^D_y \neq \emptyset$ because 
$\min(E)\leq \min(x), \min(y)$ and thus
$\hat{s}_D(x)< \min(x)$, $\hat{s}_D(y)< \min(y)$, 
$h^\rho_D(x)<\min(x)$, $h^\rho_D(y))<\min(y)$.
%$\rho_D(x)\geq \min(x)$, 
%$\min(E)\leq \rho_D(x), \rho_D(y), \max(x), \max(y)$.  
Thus  
 $h^{\rho}_D(x)=\hat{s}_D(x)$ and $h^{\rho}_D(y)=\hat{s}_D(y)$.
 Thus, trivially, $\hat{s}_D(x)=\hat{s}_D(y)< \min(E)$ 
 if and only if $h^\rho_D(x)=h^\rho_D(y)< \min(E)$.
{\bf Second}, suppose for all $x\in E^k$ of order type $ot$, 
$h^{\rho}_D(x)\geq \min(x).$
This set of order type $ot$ can be partitioned into two sets, $\{x \mid \Phi^D_x\neq \emptyset\}$ and
$\{x \mid \Phi^D_x = \emptyset\}$.  
On the first set, $\min(x) \leq h_D^\rho(x) =  \hat{s}_D(x) <\max(x)$ and on the second set $h_D^\rho(x) = \rho_D(x)\geq \min(x)$ and $\hat{s}_D(x) = \max(x) \geq \min(x)$.
Thus, $\hat{s}_D(x)\geq \min(x)$.
The same argument works if we assume for $x\in E^k$ of order type $ot$  $\hat{s}_D(x)\geq \min(x).$
Thus, for $x\in E^k$ of order type $ot$, $h^{\rho}_D(x)\geq \min(x)$ if and only if
$\hat{s}_D(x)\geq \min(x).$
\end{proof}
\end{lem}

\begin{thm}[\bfseries Regressive regularity $h^{\rho}_D$]
\label{thm:jfh}
Let $G=(N^k, \Theta)$, $r\geq 1$, $p, k\geq 2$. Let $S=\{h^{\rho}_D\mid D\subset N^k,\;|D|<\infty\}$.  Then some $f\in S$ has at most $k^k$ regressive values on some $E^k \subseteq {\rm domain}(f)=D$, $|E|=p$.
In fact, some $f\in S$ is regressively  regular over some $E$ of cardinality~$p$.
\begin{proof}
Follows from theorem~\ref{thm:jfhats} and lemma~\ref{lem:shatvsh}. 
\end{proof}
\end{thm}

It has been shown by Friedman,  Theorem~4.4 through Theorem~4.15 \cite{hf:alc}, that a special case of theorem~\ref{thm:jfh} ($\rho_D=\min$) requires the same large cardinals to prove as the jump free theorem.  
Hence, 
theorem~\ref{thm:jfh} provides a family of ZFC independent  theorems parameterized by  the $\rho_D$.

\begin{defn}[\bf $D$  capped by $E^k\subset D$]
\label{def:cap}
For $k\geq 2$, $E^k \subseteq D\subset N^k$, let $\max(D)$ be the maximum over $\max(z)$, $z\in D$.  
Let $\setmax(D)=\{z\mid z\in D, \max(z)=\max(D)\}$.
If  $\setmax(D) = \setmax(E^k)$, we say that $D$ is {\em capped by} 
$E^k\subseteq D$ with the {\em cap} defined to be $\setmax(E^k)$.
\end{defn}

The following theorem is equivalent to theorem~\ref{thm:jfh}. 
See \cite{gw:lem} for discussion and examples.

\begin{thm}[\bfseries Regressively regular $h^{\rho}_D$, capped version]
\label{thm:jfhcap}
Let $G=(N^k, \Theta)$, $r\geq 1$, $p, k\geq 2$. Let $S=\{h^{\rho}_D\mid D\subset N^k,\;|D|<\infty\}$.    Then some $f \in S$ has at most $k^k$ regressive values on some $E^k \subseteq {\rm domain}(f)=D$, $|E|=p$.
In fact, some $f \in S$ is regressively  regular over some such $E$, 
$E^k\subseteq D =  \mathrm{domain}(f)$, $D$ capped by $E^k.$
\begin{proof}
Follows from theorem~\ref{thm:jfh} by using the downward condition on $G=(N^k, \Theta)$.   
Define $D_x = \{z\mid z\in D,\,\max(z)<\max(x)\}.$
Consider $h_D^\rho$ and $E^k\subseteq D$.
Note that the downward condition on $G$  hence $G_D$ plus the recursive definition
of $h_D^\rho$ implies that $D$ can be replaced by $D_x \cup \setmax(E^k)$,
$x=\max(E)$,  without changing
the restriction $h_D^\rho \mid E^k$.
\end{proof}
\end{thm}

%\begin{defn}[\bf Diagonally restricted]
%\label{def:dagres}
%Let   $G=(N^k, \Theta),\,$, $k\geq 2,\,$ $\,F\,$ a partial selection function.
%If the set of defined values
%\[
% \{ F[z, (y_1,n_1), (y_2,n_2), \ldots, (y_r,n_r)]:\;y_i \in G^z,\, 1\leq i\leq r\}
%\]
%is empty for $z\in \diag(N^k),\,$ $r\geq 1,\,$ then $G$ is
% {\em diagonally restricted by $F$}.
%\end{defn}
%This {\em restricted diagonal} condition guarantees that for $r\geq 1$, 
%\[
%\Phi^D_z = \{ F[z, (y_1,n_1), (y_2,n_2): \ldots, (y_r,n_r)],\;y_i \in G^z_D\}
%\]
%is empty for $z\in \diag(E^k)$.
%

\begin{defn}[\bf $E$  displacement function $\delta_E$]
\label{def:cloint}
For $n\in N$, define $\gamma_E(n)$ to be the closest integer of $E=\{e_0, e_1, \ldots, e_{p-1}\}$ to $n$, ties going to the larger element of $E$.
Define $\delta_E(n) = n-\gamma_E(n)$ to be the $E$ displacement function.
\end{defn}

\begin{defn}[\bf $t$-log bounded]
\label{log:bound}
Let $p, k\geq 2$, $t\geq 1$. 
The function $\rho_D$ is 
 $t$-log bounded over $E^k\subset D$, $|E|=p,\,$ if 
$ \delta_E(\rho_D(\vec{e}_j))>0$ for $j=0, \ldots p-1$ and
$
|\{\delta_E(\rho_D(\vec{e}_j)) : \delta_E(\rho_D(\vec{e}_j))  < e_0 k^k,\, j=0, \ldots, p-1\}|
\leq t\log_2 (p).
$
We write $\rho_D\in \rm{LOG}(k,E,p,D,t)$. 

The set $\rho=\{\rho_{D}\mid \rho_D:D \rightarrow N,\, \min(x)\leq \rho_{D}(x),\,x\in D\}$
is $t$-log bounded if $\rho_D\in \rm{LOG}(k,E,p,D,t)$ when $D$ is capped by $E^k.$
We write $\rho_t$ to indicate that $\rho$ is $t$-log bounded.
\end{defn}

It is always possible to choose  $\rho_D\in \rm{LOG}(k,E,p,D,t)$. 
For example, for any $\vec{e}_j$ we can choose $\rho_D(\vec{e}_j)-e_{p-1}>0$. 
In this case, $\rho_D(\vec{e}_j)-e_{p-1} =\delta_E\rho_D(\vec{e}_j)$.
Recalling that $\rho_D(\vec{e}_j)\geq e_j$ can be arbitrarily large,
we can then choose the cardinality $|\{j:\delta_E\rho_D(\vec{e}_j) \geq e_0 k^k\}|$ large enough to make $\rho_D\in \rm{LOG}(k,E,p,D,t)$. 

\begin{thm}[\bfseries Regressive regularity of $h^{\rho_t}_D$]
\label{thm:jfhcaplog}
Let  $G=(N^k, \Theta)$, $r,t\geq 1$, $p, k\geq 2$.   
Let $S=\{h^{\rho_t}_D\mid D\subset N^k,\;|D|<\infty\}$.  
Then some $f \in S$ has at most $k^k$ regressive values on some $E^k \subseteq {\rm domain}(f)=D$, $|E|=p$.
In fact, some $f \in S$ is regressively  regular over some such $E$, 
$E^k\subseteq D =  \mathrm{domain}(f)$, $D$ capped by $E^k$ and 
$\rho_D\in \rm{LOG}(k,E,p,D,t)$.
\begin{proof}
Follows from theorem~\ref{thm:jfhcap} which states that  some $f \in S$ has at most $k^k$ regressive values on some $E^k \subseteq {\rm domain}(f)=D$, $|E|=p$.
In fact, some $f \in S$ is regressively  regular over some such $E$, 
$E^k\subseteq D =  \mathrm{domain}(f)$, $D$ capped by $E^k.$
From the definition of $\rho_t$, for each such capped pair $D$ and $E^k$,  $\rho_D$ has already been defined so that $\rho_D\in \rm{LOG}(k,E,p,D,t)$.
\end{proof}
\end{thm}

Theorem~\ref{thm:jfhcaplog} is independent of ZFC as is theorem~\ref{thm:jfhcap}.
\begin{defn}[\bf Subsets of $E$]
Define
subsets $E_l^k =\{x\mid x\in E^k, f(x)< \min(x)\}$,
$E_L^k=\{x\mid x\in E^k, f(x) < \min(E)\},\;$
%$E_U^k=\{x\mid x\in E^k, f(x)\geq \min(x)\},\;$
$\mathrm{diag}(E^k) =\{\vec{e}_0, \vec{e}_1, \ldots , \vec{e}_{p-1}\}$
where $\vec{e}_s = (e_s, \ldots, e_s).$
\end{defn}

\begin{defn}[\bf Sets of displacement]
\label{def:setsdis}
Define $S^\rho_{F,G}(k,t,E,p,D),\;$                                       %$\Delta(F, G, k, E, p, D)$, 
the {\em family of sets of displacements}, by
\begin{equation}
S^\rho_{F,G}(k,t,E,p,D),\;=\{\delta_E h_D^{\rho} E^k_l \cup \delta_E h_D^{\rho} \mathrm{diag}(E^k): k, t, F, G, E, p, D\}
\end{equation}
where $G=(N^k,\Theta),\;k\geq 2,\;$ ranges over all downward directed lattice graphs, $F:N^k\times (N^k\times N)^r\rightarrow N,\;$ $r\geq 1,\;$ ranges over all partial selection functions,
$E\subset N$ ranges over all finite subsets, $|E|=p\geq 2$, $D\subset N^k$ ranges over all  $D$ capped by $E^k,\,$ 
 $\rho=\rho_t$ is $t$-log bounded, $t\geq 1$.
%and $h_D^\rho$ has $\rho_D\in \rm{LOG}(k,E,p,D,t)$.
\end{defn}

We summarize some of the terminology.
$(1)$ $N$ the nonnegative integers. 
$(2)$ $N^k$ the nonnegative integral lattice of dimension $k\geq 2$. 
$(3)$ $\rho$ a collection of functions $\rho_D$, one for each finite $D\subset N^k$.
$(4)$ $F:N^k\times (N^k\times N)^r\rightarrow N,\;$ $r\geq 1,\;$ partial selection functions.
$(5)$ $G=(N^k,\Theta)$ a downward directed graph on $N^k$.
$(6)$ $G_D = (D, \Theta_D)$ restriction of $G$ to $D$.
$(7)$ $h_D^\rho$ functions defined recursively on $D\subset N^k$.
$(8)$ $E\subset N$, $|E|=p\geq 2$ and $\delta_E$ the $E$ displacement function.
$(9)$ $E^k\subseteq D$, $D$ capped by $E^k$.
$(10)$ $\rho_t$ a subclass of $\rho$ that are $t$-bounded defined in terms of
$\rm{LOG}(k,E,p,D,t)$.
$(11)$ $E^k_l$, $E^k_L$, $\mathrm{diag}(E^k)$ special subsets of $E^k$.

\begin{thm}[Subset sum] Let
\label{cor:subsumpoly}
$S_{F,G}^\rho(k,t, E,p,D)$ be the family of sets of displacements.
Consider the $\delta_E h_D^\rho E^k_l \cup \delta_E h_D^\rho \mathrm{diag}(E^k)
\in S_{F,G}^\rho(k,t,E,p,D)$ as instances to the subset sum problem (target 0, size measured by $p$).
For fixed $F, G, k, \rho$ consider sets of instances 
$\{\delta_E h_D^\rho E^k_l \cup \delta_E h_D^\rho \mathrm{diag}(E^k): E, p, D\}.$
For each $p$ there exists $\hat{E}$ and $\hat{D}$ such that
the subset sum problem for
$$\{\delta_{\hat{E}} h_{\hat{D}}^\rho {\hat{E}}^k_l \cup \delta_{\hat{E}}  h_{\hat{D}}^\rho \mathrm{diag}(\hat{E}^k): p=2, 3, \ldots\}$$ is solvable in time $O(p^t)$ for some $t$.

\begin{proof}
We use $\delta$ rather than $\delta_E$ to simplify the notation.
From the definition of $S_{F,G}^\rho(k,t,E,p,D)$  the parameter 
$\rho=\rho_t$ is $t$-log bounded for some $t\geq 1$.
From theorem~\ref{thm:jfhcaplog}, 
for any $p$, we  can choose ${\hat{D}}$ capped by $\hat{E}^k$ such that $h_{\hat{D}}^{\rho_t}$ is regressively regular over  $\hat{E}$.
For notational simplicity we set $\hat{E}=\{e_0, \ldots, e_{p-1}\}.$
By regressive regularity, the set $\delta h_{\hat{D}}^{\rho_t} \hat{E}^k_l$ becomes the set 
$\delta h_{\hat{D}}^{\rho_t} \hat{E}^k_L$ and,
for $x\in \hat{E}^k_L$, $\delta h_{\hat{D}}^{\rho_t}(x) = h_{\hat{D}}^{\rho_t}(x) - e_0 < 0$.
Note  $|\delta h_{\hat{D}}^{\rho_t}(x)|< e_0$ and $|\delta h_{\hat{D}}^{\rho_t} \hat{E}^k_L| < k^k$ so

\begin{equation}
\label{eq:totsum}
\sum_{x \in \hat{E}^k_L} |\delta h_{\hat{D}}^{\rho_t}(x)| < e_0 k^k.
\end{equation}

Note either
$\mathrm{diag}(\hat{E}^k) =\{\vec{e}_0, \vec{e}_1, \ldots , \vec{e}_{p-1}\}\subseteq \hat{E}_L$
or $h_{\hat{D}}^{\rho_t} \mathrm{diag}(\hat{E}^k) =\rho_{\hat{D}}\mathrm{diag}(\hat{E}^k).$ 
This follows from lemma~\ref{lem:shatvsh}, noting that for $x\in\mathrm{diag}(\hat{E}^k),\, \max(x)=\min(x)$.
In the case $\mathrm{diag}(\hat{E}^k)\subseteq \hat{E}_L,\;$ 
$\delta h_{\hat{D}}^{\rho_t} \hat{E}^k_l \cup \delta h_{\hat{D}}^{\rho_t} \mathrm{diag}(\hat{E}^k)$
consists of nonzero negative numbers and there is no solution. 
In the  case $h_{\hat{D}}^{\rho_t} \mathrm{diag}(\hat{E}^k) =\rho_{\hat{D}}\mathrm{diag}(\hat{E}^k),\,$
$\rho_{\hat{D}}\in \rm{LOG}(\hat{D},\hat{E},p,k,t)$  means (see definition~\ref{log:bound})  
\[
|\{\delta h_{\hat{D}}^{\rho_t}(\vec{e}_j) : \delta h_{\hat{D}}^{\rho_t}(\vec{e}_j)  < e_0k^k,\, j=0, \ldots, p-1,\;\}|
\leq t\log_2 (p).
\]
Thus, from equation~\ref{eq:totsum}, we can check all solutions in
$2^{k^k} 2^{t\log_2(p)} = O(p^t)$ time.
%We note that for all diagonally restricted graphs~\ref{def:dagres}, 
%$h_{\hat{D}}^\rho \mathrm{diag}(\hat{E}^k) = \rho_D\mathrm{diag}(\hat{E}^k)$
%and  the case  
%$\mathrm{diag}(\hat{E}^k) =\{\vec{e}_0, \vec{e}_1, \ldots , \vec{e}_{p-1}\}\subseteq \hat{E}_L$
%never occurs. 
\end{proof}
\end{thm}

We have proved theorem~\ref{cor:subsumpoly} from 
theorem~\ref{thm:jfhcaplog} which is independent of ZFC.
We know of no other proof.
We note that if a ZFC proof could be found that the subset sum problem is solvable in polynomial time then that result would prove theorem~\ref{cor:subsumpoly}. There would be no need for a ZFC independent proof (e.g., theorem~\ref{thm:jfhcaplog}).
In addition, if theorem~\ref{cor:subsumpoly} is itself independent of ZFC then
the polynomial time solvability of subset sum is independent of ZFC.
The intentional close relationship between theorem~\ref{thm:jfhcaplog} and 
theorem~\ref{cor:subsumpoly} leaves the possibility open that the latter is in fact independent of ZFC.

{\bf Acknowledgments:}  The author thanks  Professor Sam Buss (University of California San Diego, Department of Mathematics) and 
Professor Emeritus Rod Canfield (University of Georgia, Department of Computer Science) for numerous helpful suggestions.

%
%\cite{al:iuv}
%\cite{mt:mhy}
%\cite{sc:eth}
%\cite{jg:pos}
%\cite{hf:alc}:
%\cite{hf:nlc}
%\cite{gw:lem}

\bibliographystyle{alpha}
\bibliography{multiverse}

\end{document}